\documentclass[12pt]{article}
\usepackage{amsfonts}
\usepackage{color}
\usepackage{amssymb,amsmath,amsthm,cite,enumerate}
\usepackage{mathrsfs}

\begin{document}

\title{The $95256$-cap in $PG(12,4)$ is complete}
\author{D. Bartoli, S. Marcugini, A. Milani  and F. Pambianco \thanks{%
Partially supported by the Italian Ministero dell'Istruzione,
dell'Universit\`a e della Ricerca (MIUR) and by the Gruppo Nazionale per le
Strutture Algebriche, Geometriche e le loro Applicazioni (GNSAGA).} \\
Dipartimento di Matematica e Informatica\\
Universit\`a degli Studi di Perugia\\
Perugia (Italy)\\
\{daniele.bartoli,stefano.marcugini,alfredo.milani, \\fernanda.pambianco\}@unipg.it\\
phone\; +39(075)5855006 \\
fax~~~~~ +39(075)5855024}
\maketitle

\begin{abstract}
We describe an algorithm for testing the completeness of caps in $\mathrm{PG}(r,q)$, $%
q$ even. It allowed us to check that the $95256$-cap in $\mathrm{PG}(12,4)$ recently
found by Fu el al. (see \cite{cinesi}) is complete.
\end{abstract}

\theoremstyle{plain} \newtheorem{Theorem}{Theorem}[section] %
\newtheorem{Proposition}[Theorem]{Proposition} %
\newtheorem{Lemma}[Theorem]{Lemma} \newtheorem{Corollary}[Theorem]{Corollary}
\newtheorem{definition}[Theorem]{Definition}

{\theoremstyle{definition} \newtheorem{Definition}[Theorem]{Definition} %
\newtheorem{Example}[Theorem]{Example}}

\section{Introduction}

\bigskip

Let $\mathrm{PG}(r,q)$ be the $r$-dimensional projective space over the Galois field $\mathbb{F}_{q}$.
An $n$-cap in $\mathrm{PG}(r,q)$ is a set of points no three of which are collinear.
An $n$-cap in $\mathrm{PG}(r,q)$ is called complete if it is not contained in an $(n+1)$-cap in $\mathrm{PG}(r,q)$; see \cite{hirsh}.

The points of a complete $n$-cap in $\mathrm{PG}(r-1,q)$ can be treated as columns of a
 parity check  matrix of an  $[n,n-r,4]_{q}$ linear code with the exceptions of the complete $5$-cap in $\mathrm{PG}(3,2)$ and
the complete $11$-cap in $\mathrm{PG}(4,3)$ corresponding to the binary $[5,1,5]_{2}$ code and to the Golay $[11,6,5]_{3}$
code respectively.

An $n$-cap in $\mathrm{PG}(r,q)$ of maximal
size is called a maximal cap in $\mathrm{PG}(r,q)$. A classical problem on caps is to
determine the maximal size of complete caps in $\mathrm{PG}(r,q)$. This is also known as
the packing problem; see \cite{hirshsurvey}. Denote the size of a maximal
cap in $\mathrm{PG}(r,q)$ as $m_{2}(r,q)$, and the largest size of a known complete
cap as $\overline{m}_{2}(r,q)$.

Of particular interest is the case $q=4$, due the connection with quantum
error correction established in \cite{shor}, where a class of quantum codes,
the quantum stabilizer codes, is described in terms of certain additive
quaternary codes.

Additive quaternary codes are defined over $\mathbb{F}_{4}$ but are linear
over $\mathbb{F}_{2}$. If we restrict considering quaternary quantum codes
that are indeed $\mathbb{F}_{4}$-linear then we have the following
definition; see \cite{BierNoiSteructureQuantumCaps,BBMP2010}:
\begin {definition}
A linear quaternary quantum stabilizer code is a subspace $\mathcal{C}%
\subset \mathbb{F}_{4}^{n}$ such that $\mathcal{C}\subset \mathcal{C}^{\bot
_{H}},$ where duality is with respect to the Hermitian inner product.
\end {definition}

Here the Hermitian inner product of $x=(x_{1},...,x_{n})$ and $%
y=(y_{1},...,y_{n})$ is $\langle x,y\rangle =\sum\nolimits_{i=1}^{n}x_{i}%
\overline{y}_{i}$, where $\overline{y}=y^2$. The reason for this definition is that a linear
quaternary quantum stabilizer code $\mathcal{C}$ of length $n$, dimension $r$
and dual distance $\geq d$ (equivalently: of strength $>d$) allows the
construction of a pure quantum stabilizer code $[[n,n-2r,d]]_{4}$; see \cite[%
Theorem 1]{twistedcodes}.


A pure quantum code $[[n,n-2r,4]]$ which is linear over $\mathbb{F}_{4}$ is
obtained from a cap satisfying certain conditions; see \cite[Theorem 2.8]%
{quantgeom}:

\begin {definition}
 A cap $\mathscr{C}$ in $PG(r - 1, 4)$  is a quantum cap if it is
not contained in a hyperplane and if it satisfies the following equivalent
conditions:

\textbullet\ each hyperplane meets the cap in the same parity as the
cardinality of the cap;

\textbullet\ the corresponding quaternary $[n,r]_4$-code has all its weights
even;

\textbullet\ the corresponding quaternary $[n,r]_4$-code is self-orthogonal with respect to
the Hermitian inner product.
\end {definition}

\begin {Theorem}
 The following are equivalent:

\textbullet\ A pure stabilizer quantum code $[[n, n - 2r, 4]]$ which is linear over $\mathbb{F}_{4}$.

\textbullet\ A quantum n-cap in $PG(r - 1, 4)$.
\end {Theorem}

Much work on caps has been done, see [3-7,10-21]. The value of $m_{2}(r,4) $
is known for $k\leq 4$: $m_{2}(2,4)=6$, $m_{2}(3,4)=17$\textcolor{red}{,} and $m_{2}(4,4)=41. $

In \cite{cinesi} it is proved that $\overline{m}_{2}(8,4)=2136$, $\overline{m%
}_{2}(9,4)=5124$, $\overline{m}_{2}(10,4)=15840$, $\overline{m}%
_{2}(11,4)=36084$ and they also give a $95256$-cap in $PG(12,4)$.

Their results have been obtained by computer-supported recursive
constructions. They also present an algorithm for checking completeness of a
cap based on a bijective map between points in $PG(r,4)$ and a subset $I$ of
the positive integer set $\mathbb{N}$.

This algorithm allowed checking the completeness of the caps for $k\leq 11$,
but it is too computationally expensive for the case $k=12$. As they wrote:
\textquotedblleft But as for checking completeness of larger caps in $%
PG(r,4)$, $r\geq 12$, new algorithms are needed.\textquotedblright ; see
\cite[Section 5]{cinesi}. We propose a new fast algorithm that allowed to
face also this case: we verified that the $95256$-cap in $PG(12,4)$ is
complete, so $\overline{m}_{2}(12,4)=95256$. Our algorithm is based on a
compact representation of the points of $\mathrm{PG}(r,q)$, $q$ even\textcolor{red}{,} and on
minimizing the computational costs of the operations more often performed
during the check of the completeness of the cap.

Section 2 describes the algorithm and applies it in $PG(12,4)$. Section 3
contains the generalization of the algorithm to other even values of $q$ and
other dimensions.

\section{A new algorithm for checking completeness of a cap}

In \cite{cinesi} an algorithm for checking completeness of a cap $\mathscr{C}
$ in $PG(r,4)$ is presented. It is based on a bijective map $\phi $\ between
points in $PG(r,4)$ and a subset $T(r)$ of the positive integer set $\mathbb{%
N}$:

\begin{center}
$\phi :PG(r,4)\rightarrow T(r)$,

$\phi :P\mapsto \phi (P)$,
\end{center}

where

\begin{center}
$P=(x_{0},x_{1},...,x_{r})^{T}$,

$\phi (P)=4^{r}x_{0}+4^{r-1}x_{1}+$\textperiodcentered \textperiodcentered
\textperiodcentered $+4x^{r-1}+x^{r}$.
\end{center}

It can be easily seen that a cap is complete if and only if each point $P$ of $\mathrm{PG}(r,q)$ not belonging to the cap lies on a
secant line of the cap. In this case we say that $P$ is covered.

To keep track of the covering of the points, a vector $U$ of size $|T(r)|$ is
used. Initially all elements of $U$ are set to be 1.

Then all pairs of points of the cap are considered. For each pair of points  $(P_{i}, P_{j})$ the three other points belonging to the line through $%
P_{i}$, $P_{j}$ are computed. To do this all linear combination $Q=\alpha
P_{1}+P_{2},\alpha \in \mathbb{F}_{4}\backslash \{0\}$ are computed. The
point $Q$ is normalized, choosing a representation with the leftmost non-zero
coordinate equal to $1$. Finally the position $\phi (Q)$\ of $U$\ is set to $%
0$. The process continues until all elements of $U$ became $0$ or all pairs
of points of the cap have been considered.

At the end the cap is complete if and only if all elements of $U$ are $0.$
Table 1 reports the time cost of the algorithm using an Intel(R) Xeon(R) CPU
E5504 @ 2.00GHz; see \cite[Table 1]{cinesi}. However, the paper does not mention the unit of time used in Table 1.

\begin{table}[tbp]
\caption{The time and space cost of the algorithm of \protect\cite{cinesi}}
\begin{center}
\begin{tabular}{|c||c|c|c|c|}
\hline
Size of cap & 2136 & 5124 & 15840 & 36084 \\ \hline
Time & 7805 & 46244 & 428029 & 2261301 \\ \hline
\end{tabular}%
\end{center}
\end{table}

We devised a new algorithm for checking the completeness of a cap in $PG(r,4)$
choosing a representation that optimizes the computational cost of the main
operations of the previous algorithm: the computation of $Q=\alpha
P_{1}+P_{2} $ and the normalization of a point $Q$.

Let be $\mathbb{F}_{4}=\ \{0,1,\omega ,\overline{\omega }\}$, where $\omega^2=\overline{\omega }$, $\overline{\omega }=\omega+1$\textcolor{red}{,} and $\omega^3=1$. If we define a
representation function $\rho :$ $\mathbb{F}_{4}\rightarrow \mathbb{N}$ as in
the following:
$$\rho (0)\rightarrow 0, \; \rho (1)\rightarrow 1, \; \rho (\omega )=2, \; \rho (\overline{\omega })=3,$$
then we have
$$a+b=\rho (a) \symbol{94} \rho(b), \;\; a,b \in \mathbb{F}_{4},$$
  where $\symbol{94}$ is the bitwise exclusive or operator.

Moreover\textcolor{red}{,} if $P=(x_{0},x_{1},...,x_{r})^{T}$ then the binary representation
of $\phi (P)$\ is $\rho (x_{0})\rho (x_{1})\ldots \rho (x_{r})$.

It means that if $Q=P_{1}+P_{2}$ then $\phi (Q)=\phi (P_{1})\symbol{94}\phi
(P_{2})$.

This allows the computation of the sum of two points of $PG(r,4)$ by one
integer operation.

The multiplication of one point $P$ by a scalar is applied only to the
points of the cap. It can be pre-computed before the beginning of the
check for completeness, so at the cost of having a data structure of size $3|%
\mathscr{C}|$ all multiplications  by a scalar are avoided.

The other expensive operation of the algorithm in \cite{cinesi} is the normalization of a point $P$. It should be to computed each
time $Q=\alpha P_{1}+P_{2}$, $P_{1},P_{2}\in \mathscr{C}$ is computed, i.e. $%
3/2|\mathscr{C}|^{2}$ times. We propose a trade-off between computational time and
memory space: we use a vector $U$ of size $3|PG(r,4)|$ to keep trace of
the fact that a point $Q$ is covered by $\mathscr{C}$ or not; initially all
elements of $U$ are set equal to $0$. When $Q=\alpha P_{1}+P_{2}$ is computed,
then the element $\phi (Q)$ is set equal to $1$ without before normalizing $Q$.
In this way all the $3/2|\mathscr{C}|^{2}$ normalization operations are avoided.
At the end, when the covering of all points of $PG(r,4)$ is tested, first
the normalized form of a point $Q$ is tested checking the element of $U$ of
position $\phi (Q)$; if it is not covered also $\phi (\omega Q)$ and $\phi (%
\overline{\omega }Q)$ are checked: if any of
$\phi (Q), \phi (\omega Q), \phi (\overline{\omega }Q)$ is equal to 1, then $Q$ is  covered.

Let be $n=|\mathscr{C}|$, $m=|PG(r,q)|$, $i=$ the size of an integer, $c=$
the size of a character. The total cost of our algorithm is:
\begin{center}
space: $3n\cdot i+3m\cdot c+c_{1};$\\
time: $3c_{2}n+3/2c_{2}n^{2}+3c_{3}m\cdot c+c_{4}$;
\end{center}

 \noindent where $c_1,\ldots ,c_{4}$ are constants.

The algorithm has been implemented in C language.

Table 2 reports the time and space cost of the algorithm using an Intel(R)
Core(TM) i7-4510U CPU @ 2.00GHz; space is measured in megabyte, while time
is measured in milliseconds.

We tested the completeness of the caps presented in \cite{cinesi} in $%
PG(r,q) $, $r = 8, \dots, 12$. Note that for constructing the $5124$-cap in $%
PG(9, 4) $ and the $36084$-cap in $PG(11, 4)$ we were not able to obtain a
cap following the selection of columns suggested in \cite[Section 3.2]%
{cinesi}. For the $5124$-cap we had to exclude the vectors for $j \in \{
2,14,24,25\}$ instead of $j \in \{2,14,15,24\}$ as suggested in the paper,
whereas for the $36084$-cap we had to exclude the vectors for $i \in
\{1,\dots,16,271\}$ instead of $i \in \{1,256,\dots,271\}$ as suggested in
the paper.
\begin{table}[tbp]
\caption{The time and space cost of the new algorithm}
\begin{center}
\begin{tabular}{|c||c|c|c|c|c|}
\hline
Size of cap & 2136 & 5124 & 15840 & 36084 & 95256 \\ \hline
Space (Mb)& 5 & 18 & 76 & 324 & 1382 \\ \hline
Time (Milliseconds)& 32 & 78 & 707 & 3574 & 98321 \\ \hline
\end{tabular}%
\end{center}
\end{table}

\section{Generalization of the algorithm}

In the previous section we applied our algorithm in $PG(12,4)$. The key
idea, choosing a representation for the elements of $\mathbb{F}_{q}$ and the
points of $\mathrm{PG}(r,q)$\ that minimize the computational cost of the operations
most often performed during the test of completeness of a cap, can be
applied for every even $q$.

When testing the completeness of a cap $\mathscr{C}$ the value $Q=\alpha
P_{1}+P_{2}$, $P_{1},P_{2}\in \mathscr{C}$, $\alpha \in \mathbb{F}%
_{q}\backslash \{0\}$ has to be computed.

To avoid to compute the same value $\alpha P$ several times, it is convenient
to compute it at the beginning of the algorithm and store the results.
Therefore the main operation to compute remains the sum between two vectors
representing points of $\mathrm{PG}(r,q)$.

When considering a representation of $\mathbb{F}_{2^{k}}$, we can either choose
a form that facilitate the computation of multiplication (we see the
non-zero elements of $\mathbb{F}_{2^{k}}$ as powers of the primitive element) or
can choose a form that facilitate the computation of addition (we see the
the elements of $\mathbb{F}_{q}$ as polynomials of $\mathbb{F}_{2}[X]$ of degree less than $%
k $; addition is defined in the natural way, whereas multiplication is defined
modulo a fixed irreducible polynomial of degree $k$).

We choose the latter representation and define $\rho :$ $\mathbb{F}_{2^{k}}\rightarrow \{0, \dots, 2^{k}-1 \}$ as $\rho
:p(x)\mapsto p(2)$.
We have that $\rho (p(x)+q(x))=\rho (p(x))\symbol{94}\rho (q(x))$, where $%
\symbol{94}$ is the bitwise exclusive or. It means that in this
representation addition on $\mathbb{F}_{2^{k}}$\ reduces to one bitwise
arithmetic operation on integers.

Moreover we can define a representation of the points of $PG(r,2^{k})$ in
the following way:

\begin{center}
$\phi :PG(r,2^{k})\rightarrow  \mathbb{N}$,\\
$\phi :P\mapsto \phi (P)$,
\end{center}

where

\begin{center}
$P=(x_{0},x_{1},...,x_{r})^{T}$,

$\phi (P)=(2^{k})^{r}\rho (x_{0})+(2^{k})^{r-1}\rho (x_{1})+\dots +(2^{k})\rho (x^{r-1})+\rho (x^{r})$\textcolor{red}{.}
\end{center}

In this way a point $P$ of $PG(r,2^{k})$\ is represented by an integer $n$.
If we consider the binary representation of $n$, the coordinate $x_{i}$\ is
represented by the bits of $n$ in position $(r-i)k+1\ldots (r-i+1)k$ that
are the binary representation $\rho (x_{i})$. To compute $\phi (Q)$, $%
Q=P_{1}+P_{2}$ it is sufficient computing $\phi (P_{1})\symbol{94}\phi
(P_{2})$, just one bitwise arithmetic operation. In a real implementation,
usually an (unsigned) integer has a $32$ bit representation, so $\phi (P)$
can be represented by a single integer if $kr\leq 32$, otherwise more
integers are needed.

Our algorithm trades computational time for memory space. Let be $n=|%
\mathscr{C}|$, $m=|PG(r,2^{k})|$, we need $n(2^{k}-1)$ integers
to represent $\alpha P$, $P\in \mathscr{C}$, $\alpha \in \mathbb{F}%
_{2^{k}}\backslash \{0\}$ and $m(2^{k}-1)$ booleans to represent the fact
that $\alpha P$, $P\in PG(r,2^{k})$, $\alpha \in \mathbb{F}%
_{2^{k}}\backslash \{0\}$ is saturated or not. As $n<m$, the latter value is
more relevant. If the memory requested by the algorithm is too big, then
memory space can be traded for computational time. For example the set of
points of $PG(r,2^{k})$ can be divided into $s$ subsets small enough to be
represented. Than the test for completeness can be repeated $s$ times; each
time the completeness of the points of one subset is tested. Note that the
computations for the different subsets are independent, so they can be
performed in parallel. This is a form of parallelism based on the splitting
of data: it is simple and effective.

\section{Conclusion}

We presented an algorithm for checking the completeness of a cap $\mathscr{C}
$ in $PG(r,2^{k})$. The key idea is making the more effective as possible
the operations that are performed more often.

We identified the following operations:

\begin{enumerate}
\item computing the points belonging to a line passing through two points of $\mathscr{C}$;
\item normalizing a point $P$ finding a representation with the leftmost
non-zero coordinate equal to 0;
\item keeping track if a point is saturated or not.
\end{enumerate}

1. is obtained representing the points of $PG(r,2^{k})$ as integers, so that
the sum of two points is computed as a bitwise exclusive or between
integers. Multiplication by a scalar for the points of $\mathscr{C}$\ is
precomputed at the beginning of the algorithm.

2. is avoided by the way 3. is performed; we use a vector $V$ of booleans
that represents each point $P$ of $PG(r,2^{k})$ $2^{k}-1$ times: we
represent all vectors $\alpha P$, $\alpha \in \mathbb{F}_{2^{k}}\backslash
\{0\}$; $P$ is saturated is at least one of the $\alpha P$\ is saturated.

We used this algorithm for proving that the $95256$-cap in $PG(12,4)$ of
\cite{cinesi}\ is complete.

When using this algorithm with greater values of $r$ or $k$, it can became
too memory consuming. In this case one either can reintroduce point 2. and
represent in point 3. only the normalized version of $P$ or can split $V$ in
subvectors $V_{i}$ and test the covering of the points of $V_{i}$
separately. In the latter way memory is saved, but computational time increases.
However these computations can be performed in parallel.

\end{document}